\documentclass{article}

\usepackage{amsfonts}
\usepackage{amssymb}
\usepackage{amstext}

\usepackage{theorem}
 \newtheorem{thm}{Theorem}[section]
 \newtheorem{cor}[thm]{Corollary}
 \newtheorem{lem}[thm]{Lemma}
 \newtheorem{prop}[thm]{Proposition}
 \theoremstyle{definition}
 \newtheorem{defn}[thm]{Definition}
 \theoremstyle{remark}
 \newtheorem{rem}[thm]{Remark}
 \newtheorem{ex}{Example}

\newenvironment{proof}
	{\par {\it Proof:}}
 	{\hfill $\square$ \medskip}

\newcommand\ds{\displaystyle}
\renewcommand\l{\lambda}
\renewcommand\t{\times}
\newcommand\<{\langle}
\renewcommand\>{\rangle}
\newcommand\A{{\mathcal A}}
\newcommand\B{{\mathcal B}}
\renewcommand\H{{\mathcal H}}
\newcommand\V{V}
\newcommand{\wlim}{\mathop{\mbox{\rm w-lim}}}

\begin{document}
\title{Eigenfunction Expansions and \\Transformation Theory}

\author{Manuel Gadella$^1$ and Fernando G\'omez$^2$}

\maketitle

\begin{center}
\small
$^1$Dpto. de F\'{\i}sica Te\'orica. Universidad de Valladolid.\\
Facultad de Ciencias, Prado de la Magdalena, s.n.,\\
47005 Valladolid, Spain.\\
email: {\tt gadella@fta.uva.es}\\
$^2$Dpto. de An\'alisis Matem\'atico. Universidad de Valladolid.\\
Facultad de Ciencias, Prado de la Magdalena, s.n.,\\
47005 Valladolid, Spain.\\
email: {\tt fgcubill@am.uva.es}
\end{center}

\begin{abstract}
Generalized eigenfunctions may be regarded as vectors of a basis in a particular  direct integral of Hil\-bert spaces or as elements of the antidual space $\Phi^\times$ in a convenient Gelfand triplet $\Phi\subseteq{\mathcal\H}\subseteq\Phi^\times$. 
This work presents a fit treatment for computational purposes of transformations formulas relating different generalized bases of eigenfunctions in both frameworks direct integrals and Gelfand triplets. 
Transformation formulas look like usual in Physics li\-te\-ra\-tu\-re, as limits of integral functionals but with well defined kernels. Several approaches are feasible. Here Vitali and martingale approaches are showed.
\end{abstract}

{\small


}

%

\section{Introduction}

Eigenfunction expansions appear in the most varied domains as for example in the basis of Dirac formulation of Quantum Mechanics \cite{D}, where  each complete set of
commuting observables (csco) $A_1,A_2,\dots,A_n$ is supposed to have a generalized basis of kets
$|\l_1,\l_2,\dots,\l_n\>$ satisfying the
following properties:

\smallskip

1.- The kets $|\l_1,\l_2,\dots,\l_n\rangle$ are
generalized eigenvectors of the observables $A_1,A_2,\dots,A_n$,
i.e.,
$$
A_j\,|\l_1,\l_2,\dots,\l_n\>=\l_j\,|\l_1,\l_2,\dots,\l_n\>\,,\quad
(j=1,2,\dots,n)\,,
$$
where $\l_j$ is one of the possible outcomes of a measurement
of the observable $A_j$, $j=1,2,\dots,n$. Let us call $\Lambda_j$
the set of these outcomes.

\smallskip

2.- For each pure state $\varphi$, one has the following integral
decomposition:
$$
\varphi=\int_\Lambda
\<\l_1,\l_2,\dots,\l_n|\varphi\>\,|\l_1,\l_2,\dots,\l_n\>\,
d\l_1d\l_2\cdots d\l_n\,,
$$
where $\Lambda=\Lambda_1\times\Lambda_2\times\dots\times\Lambda_n$
is the Cartesian product of the $\Lambda_j$ and for (almost) all
$(\l_1,\l_2,\dots,\l_n)\in\Lambda$,
$\<\l_1,\l_2,\dots,\l_n|\varphi\>$ is a
complex number that gives the coordinates of $\varphi$ in the
generalized basis $|\l_1,\l_2,\dots,\l_n\>$.

\smallskip

Following von Neumann \cite{VN}, observables in Quantum Mechanics
are represented by selfadjoint operators in a Hilbert space $\H$. 
Therefore, a csco $A_1,\dots,A_n$ is given by a set of $n$
self adjoint operators, also called $A_i$, whose respective Hilbert space spectra are
the sets $\Lambda_i$. The
classical version of the spectral theorem asso\-cia\-tes to the family
of selfadjoint operators $A_i$ a family of
commuting Borel spectral measure spaces $(\Lambda_i,\B_i,\H,P_i)$, where $i=1\ldots,n$; see Section \ref{app1}. 

Consider $\Lambda$ as a metric space with the product
topology and the co\-rres\-ponding Borel
$\sigma$-algebra $\B_\Lambda$. Since the
commuting spectral measures are Borel, there exists a unique
spectral measure space of the form $(\Lambda,\B_\l,\H,P)$ such
that
$$
P(\Lambda_1\times\ldots\times
E_i\times\ldots\times\Lambda_n)=P_i(E_i),\quad (E_i\in\B_i,\,\,i=1\ldots,n).
$$
The spectral measure $P$ is called the {\it product} of $P_1,\ldots,P_m$ \cite{BS}. 
Thus, giving a csco $\{A_1,A_2,\dots,A_n\}$ is equivalent
to give the corresponding Borel spectral measure space $(\Lambda,\B_\Lambda,\H,P)$.

It is well known that only for the pure discrete part of the spectrum there exist eigenvectors belonging to $\H$. 
For the continuous part, generalized eigenvectors or eigenfunctions have been considered as 
components of an orthonormal measurable basis in a profitable direct integral decomposition 
of the Hilbert space $\H$ \cite{VN49,MA} or as elements of the space $\Phi^\t$ of a conveniently chosen rigging 
\begin{equation}\label{0}
\Phi\subseteq H\subseteq\Phi^\t,
\end{equation}
where $\Phi$ is a dense subspace of $\H$ with its own topology $\tau_\Phi$
and $\Phi^\t$ is its topological (anti)dual space,  i.e. the space of
$\tau_\Phi$-continuous antilinear forms on $\Phi$. A triplet of the form (\ref{0}) is usually called
a {\it Gelfand triplet} or {\it rigged Hilbert space} (RHS).
Gelfand \cite{GK,GS,GV} was the first to give a precise meaning to generalized eigenvectors, which was later elaborated, among others, by Berezanskii \cite{BERE65,BSU}, Maurin \cite{MAUR68}, Foia\c s \cite{F}, Roberts \cite{R}, Melsheimer \cite{ME} and Antoine \cite{A,AI}. For a detailed exposition see previous works by the authors \cite{GG,GG1}.
Such mathematics are being useful on a great variety of areas in Physics. In particular, RHS have been used in the study of Scattering Theory \cite{KK70}, resonances \cite{B,BG,AI,AGP,AP,BOIV,CBPR}, singular states in Statistical Mechanics \cite{VH,AS,CGIL}, spectral decompositions associated to chaotic maps \cite{AT,SATB}, time irreversibility \cite{PPT,AP} or axiomatic theory of quantum fields
\cite{BLT}.

The aim of this paper is to discuss a particular
aspect of eigenfunction ex\-pan\-sion theory: the transformations formulas relating different generalized bases. 
Given two different csco 
$A_1,A_2,\dots,A_n$ and $B_1,B_2,\dots,B_n$ with respective ge\-ne\-ra\-li\-zed bases $|\l_1,\l_2,\dots,\l_n\>$ and $|\xi_1,\xi_2,\dots,\xi_n\>$, Dirac \cite{D} introduced transformation formulas relating coordinates of a pure state $\varphi$ in both bases as follows:
\begin{equation}\label{1g}
\begin{array}{l}
\<\l_1,\l_2,\dots,\l_n|\varphi\>= 
\\
\ds =\int \<\l_1,\l_2,\dots,\l_n|\xi_1,\xi_2,\dots,\xi_n\>\,
\<\xi_1,\xi_2,\dots,\xi_n|\varphi\>\,d\xi_1 d\xi_2\cdots d\xi_n\,,  
\end{array}
\end{equation}
\begin{equation}\label{2g}
\begin{array}{l}
\<\xi_1,\xi_2,\dots,\xi_n|\varphi\>= 
\\
\ds=\int \<\xi_1,\xi_2,\dots,\xi_n|\l_1,\l_2,\dots,\l_n\>\,
\<\l_1,\l_2,\dots,\l_n|\varphi\>\, d\l_1 d\l_2\cdots d\l_n \,. 
\end{array}
\end{equation}
These linear formulas are straightforward ge\-ne\-ra\-li\-za\-tions of the familiar ones on finite-dimensional Hilbert spaces. In infinite-dimensional Hilbert spaces, however, the situation is not so simple, because the integral kernels 
$$
\<\l_1,\l_2,\dots,\l_n|\xi_1,\xi_2,\dots,\xi_n\>\quad \text{ and }\quad
\<\xi_1,\xi_2,\dots,\xi_n|\l_1,\l_2,\dots,\l_n\>
$$ 
have in general no mathematical meaning.
A preliminary discussion of the
problem is given in the early literature mentioned above, in particular, the second part of Melsheimer \cite{ME}. 
This work intends giving to transformation theory a new mathematical treatment. A fit treatment for computational purposes that can be described as follows: 
For two spectral measure spaces on the same separable Hilbert space, 
$(\Lambda,{\A},\H,P)$ and $(\Xi,{\B},\H,Q)$,
transformation formulas between the respective eigenfunction expansions are obtained as limits of integrals like those of Equations (\ref{1g}) and (\ref{2g}) but with well defined kernels.
These approximate integral functionals are provided by Vitali approach to Radon-Nikodym derivatives (see Appendix), which, in particular, is sufficient to deal with the absolutely continuous part of the spectral decomposition of a cso. 
This is done in both direct integral and rigged Hilbert space frameworks, where Vitali approach leads to the functionals $\tilde\Gamma^\t_{nl}$ and $\Gamma^\t_{nl}$ given in Equations (\ref{fun.ant.Gammat.deg.nc}) and (\ref{fun.ant.Gammaxx}), respectively.
In both frameworks pointwise or weak convergence of approximate functionals
is assured (Theorems \ref{p18} and \ref{p16}).
Convergence with respect to finer topologies is only possible in rigged Hilbert spaces:
uniform convergence on precompact sets is proved for barrelled tvs $\Phi$ (Corollary \ref{p5}), whereas convergence with respect to strong topology $\beta(\Phi^\t,\Phi)$ is verified on barrelled and
semi-reflexive tvs $\Phi$, in particular when $\Phi$ is a Montel
space, a Montel-Fr\'echet space or a tvs with nuclear strong dual
(Corollary \ref{p9}).
The paper is organized as follows: Section \ref{sc2} collects some notions and results 
used along the work relative to spectral measure spaces, direct integrals of
Hilbert  spaces, locally convex equipments and eigenfunction expansions.
Section \ref{sctf} contains Vitali approach to transformation theory described above. Some comments about other approaches using martingale theory or generalized Cauchy-Stieltjes and Poisson integrals are included in Section \ref{sect.32}. Finally, an Appendix added at the end reviews briefly Vitali systems. 


\section{Preliminaries. Eigenfunction Expansions.}\label{sc2}

In this Section we introduce the terminology and preliminary results 
related with the mathematical structures used along this paper: spectral
measure spaces, direct integrals of Hilbert spaces and locally convex equipments of
spectral measures.  For the first two structures the terminology is mainly that of 
Birman-Solomjak \cite{BS}. For equipments and eigenfunction expansions see previous works by the authors \cite{GG, GG1}.

\subsection{Spectral Measure Spaces}\label{app1}

Let $(\Lambda,{\A})$ be a measurable space, let $\H$ be a separable 
Hilbert space with scalar product $(\cdot,\cdot)$ and norm $||\cdot||$ (we  
consider separable Hilbert spaces only) and let ${\mathcal P}={\mathcal P}(\H)$ be the set 
of orthonormal projections on $\H$. A {\it spectral measure} on $\H$ is a mapping
$P:{\A}\to {\mathcal P}$ satisfying the following conditions:
(1) {\it countable additivity}, i.e. if $\{E_n\}$ is a finite or countable
set  of disjoint sets of ${\A}$, then $P(\cup_n E_n)=\sum_n P(E_n)$ in
strong sense;
(2) {\it completeness}, i.e. $P(\Lambda)=I$.
Then $(\Lambda,{\A},\H,P)$ is referred as a {\it spectral measure space}.

If $\Lambda$ is a complete separable metric space and  ${\A}$ is the
$\sigma$-algebra of all Borel subset of $\Lambda$, then $P$ is called a
{\it Borel spectral  measure}.
The classical version of the spectral theorem \cite{VN} associates to 
each normal operator $A$ defined on $\H$ a Borel spectral measure space
$(\Lambda,{\A},\H,P)$, being $\Lambda=\sigma(A)$, the spectrum of
$A$, and  $P(E)$ the orthogonal projection corresponding to $E\in\mathcal
A$.

Here we consider the more general situation introduced in \cite{BS}: measure spaces
$(\Lambda,{\A},\mu)$ such that $\mu$ is a $\sigma$-finite measure with a countable basis 
(a sequence $\{E_n\}$ of measurable sets
in $\A$ such that for any $E\in\A$ and for all $\varepsilon>0$,
there exists an $E_k$, from the sequence $\{E_n\}$, such that $\mu[(E\backslash E_k)\cup
(E_k\backslash E)]<\varepsilon$).

By the {\it type} $[\mu]$ of a measure $\mu$ defined on $(\Lambda,\A)$
we understand the equivalence class of all measures $\nu$ on $(\Lambda,\A)$ such that $\mu$ and $\nu$ are mutually absolutely continuous, that is, they have the same class of null sets.
If the measure $\mu$ is absolutely continuous with
respect to the measure $\nu$, we write $[\nu]\succ [\mu]$.

Given a spectral measure space $(\Lambda,{\A},\H,P)$,  every
pair $f,g\in\H$ define the complex measure: 

$$
\mu_{f,g}(E):=(f,P(E)g),\quad (E\in \A).
$$
When $f=g$ we write $\mu_f:=\mu_{f,f}$ and therefore,
$\mu_f(E)=(f,P(E)f)$. 
We say that a nonzero vector $g\in\H$ is of {\it maximal type} with respect to the
spectral measure $P$ if for each $f\in\H$,  $[\mu_g]\succ [\mu_f]$. 
In this case, $g$ is called a maximal vector. Such maximal vectors 
always exist provided that $\H$ is separable. 
The type $[\mu_g]$ of a maximal vector is called the  {\it spectral
type} of $P$ and denoted by $[P]$.

For an element $g$ of $\H$ we denote by $\H_g$ the closed subspace 
of $\H$ defined by 
$$
\H_g={\rm adh\,}\{f\in\H: f=P(E)g\},
$$ 
where $E$ runs over $\A$ and {\it adh} denotes adherence.
A family of vectors $\{g_j\}_{j=1}^m$ in $\H$, where $m\in
\{1,2,\dots,\infty\}$,
is called a {\it generating system} of $\H$ with respect to $P$ if  
$\H$ is the orthogonal sum of the spaces $\H_{g_j}$:
\begin{equation}\label{odh}
\H=\bigoplus_{j=1}^m\H_{g_j}.
\end{equation}
There exists a generating system $\{g_j\}_{j=1}^m$ such that
\begin{equation}\label{seq.ty}
[P]=[\mu_{g_1}] \succ [\mu_{g_2}] \succ\dots
\end{equation}
For each $k\in\{1,2,\ldots,m\}$ let $\Lambda(g_k)$ be the support of $\mu_{g_k}$.
The family $\{\Lambda_k\}_{k=1}^m$ of subsets of $\Lambda$ given by
$$
\begin{array}{rcl}
\Lambda_k &=&\Lambda (g_k)\backslash \Lambda  (g_{k+1})\,, \hskip0.6cm
(k<m),
\\[1ex]
\Lambda_m &= &\bigcap_{k\in \{1,2,\ldots,m\}} \Lambda(g_k), 
\end{array}
$$
is a partition of $\Lambda$ and then we can define the  {\it multiplicity
function} $N_p$ of $P$ as
$$
N_P(\l)=k\,,\hskip0.6cm {\rm if}\hskip0.6cm \l\in\Lambda_k
$$
The sequence of types (\ref{seq.ty}) and the multiplicity function $N_p$ are the unitary invariants of $P$.

Now, let $S(\Lambda,P)$ be the set of complex-valued measurable functions on $\Lambda$ 
which are finite almost
everywhere with respect to the measures $\mu_g(E)=(g,P(E)g)$ for all 
$g\in\H$ or, equivalently, with respect to $\mu_g$ for some (all) $g\in\H$ of maximal type (we say that
these functions are finite almost everywhere with respect to $P$  or
that are
$[P]$-a.e. finite). For each $\phi\in S(\Lambda,P)$, we can define the
operator

$$
J_\phi:= \int_\Lambda \phi\,dP 
$$
with domain
${\mathcal D}_\phi=\{f\in\H : \int_\Lambda |\phi|^2\,d\mu_f<\infty\}$.
The main properties of the operators $J_\phi$ are listed in the 
following \cite{BS}:

\begin{prop}\label{pA1}
If $f,g\in\H$ and $\phi\in S(\Lambda,P)$, then: 

\smallskip
(i) $(f,J_\phi g)=\int_\Lambda \phi\,d\mu_{f,g}$, if $g\in {\mathcal D}_\phi$.

\smallskip
(ii) $(f,J_\phi f)=\int_\Lambda \phi\,d\mu_f$, if $f\in {\mathcal D}_\phi$.

\smallskip
(iii) $||J_\phi f||^2=\int_\Lambda |\phi|^2\,d\mu_f$, if $f\in  {\mathcal
D}_\phi$.

\smallskip
(iv) $||J_\phi||= [P]\mbox{-}{\rm ess\, sup\,} |\phi|$, if  $\phi\in
L^\infty(\Lambda,[P])$.

\smallskip
(v) $g\in{\mathcal D}_\phi$ if and only if $\phi\in L^2(\Lambda,\mu_g)$.

\smallskip
(vi) $\H_g$ reduces $J_\phi$, i.e.,
$J_\phi(\H_g\cap {\mathcal D}_\phi)\subset \H_g$ and
$J_\phi(\H_g^\bot\cap {\mathcal D}_\phi)\subset \H_g^\bot$.

\smallskip
(vii) The mapping 
$V_g: L^2(\Lambda,\mu_g)\rightarrow \H_g :\phi \mapsto J_\phi g$
is unitary.

\smallskip
(viii) 
$\H_g= \{ J_\phi g:\phi\in S(\Lambda,P),g\in{\mathcal D}_\phi\}$.
\end{prop}


\subsection{Direct Integrals of Hilbert Spaces}\label{app2}

Let $(\Lambda,{\A})$ be a measurable space. We say that a family 
$\{{\H}_\l\}_{\l\in \Lambda}$ of separable Hilbert spaces, with
scalar product $(\cdot,\cdot)_\l$ and norm $||\cdot||_\l$, together with
a countable set $\{u_j(\l)\}_{j=1}^\infty$ of elements of the product
$\prod_{\l\in \Lambda} {\H}_\l$ is a {\it measurable family of Hilbert
spaces on $(\Lambda,{\A})$} if for every $j,k\in{\mathbb N}$ the functions 
$( u_j(\l),u_k(\l))_{\l}$ are measurable and for each $\l\in \Lambda$ 
the subspace
${\rm span}\{u_j(\l):1\leq j<\infty\}$ is dense in ${\H}_\l$.
In such case the {\it function of dimension}
$N(\l):={\rm dim}({\H}_\l)$ is measurable and we say that an element 
$f$ of $\prod_{\l\in \Lambda} {\H}_\l$ is a {\it measurable vector
family} if the function 
$( f(\l),u_j(\l))_{\l}$ is measurable for each $j$.
Given a measure $\mu$ on $(\Lambda,{\A})$, the space of measurable vector families $f$ on ${\A}$ such that 
$\int_\Lambda ||f(\l)||^2_{\l}\,d\mu(\l)<\infty$ is called the  {\it direct integral}
of the spaces ${\H}_\l$ with respect to $\mu$ and denoted by 
$$
{\H}_{\mu,N}\quad \text{or}\quad \int_\Lambda^\oplus {\H}_\l\,d\mu(\l).
$$
If, as usual, we identify functions that coincide $\mu$-a.e., the space 
${\H}_{\mu,N}$ is a Hilbert space with the inner product
$(f,g)_{\H_{\mu,N}}:=\int_\Lambda \big( f(\l),g(\l)\big)_{\l}\,d\mu(\l)$.
This structure does not depend on $\{u_j\}$.
The sets
$\Lambda_k=\{\l\in \Lambda:N(\l)=k\}$, where $k=1,2,\ldots,\infty$,
are measurable and there exist sequences $\{e_k(\l)\}_{k=1}^\infty$ of 
measurable vector families with the following properties:

\smallskip

(a) for each $\l$, $\{e_k(\l)\}_1^{N(\l)}$ is an orthonormal basis  for
${\H}_\l$, and $e_k(\l)=0$ for $k>N(\l)$;

\smallskip

(b) for each $k$ there is a measurable partition of $\Lambda$, 
$\Lambda=\bigcup_{l=1}^\infty \Lambda_l^k$, such that on each
$\Lambda_l^k$, $e_k(\l)$ is a finite linear combination of the $u_j(\l)$
with coefficients depending measurably on $\l$.

\smallskip\noindent
Such a sequence $\{e_k\}$ is called a {\it orthonormal measurable basis} 
of the direct integral ${\H}_{\mu,N}$ and its choice determines a
unitary isomorphism between ${\H}_{\mu,N}$ and 
$L^2(\Lambda_\infty,\mu;l^2)\oplus\big(\oplus_1^\infty
L^2(\Lambda_m,\mu;{\mathbb C}^m)\big)$.

Any measurable set $ E \in{\A}$ generates the operator  $\hat{P}(E)$
on ${\H}_{\mu,N}$, multiplication by the characteristic function
$\chi_ E $ of $ E $ ($\chi_E(\l)$ is one if
$\l\in E$ and zero otherwise),

$$ 
\hat{P}( E )g:=\chi_ E  g,\quad (g\in{\H}_{\mu,N}).
$$
The family of projections $ \hat{P}( E )$, ($E
\in{\A}$), defines a spectral measure on ${\H}_{\mu,N}$. 
The equalities $\mu(E)=0$ and $\hat{P}(E)=0$ are equivalent and,
thus, the {\it types} of $\mu$ and $\hat{P}$ coincide and also
$S(\Lambda,\mu)=S(\Lambda, \hat{P})$ and
$L^\infty(\Lambda,\mu)=L^\infty(\Lambda,\hat{P})$. 
In general, every $\phi\in S(\Lambda,\mu)$ defines a 
{\it multiplication or diagonal operator} 
on ${\H}_{\mu,N}$, which we denote by $Q_\phi$, given by

$$
\begin{array}{l}
(Q_\phi f)(\l):=\phi(\l)f(\l),\quad \forall f\in D(Q_\phi),\\[1ex]
D(Q_\phi)=\{f\in {\H}_{\mu,N}:\int|\phi(\l)|^2||f(\l)||^2_{\l}\,
d\mu(\l)<\infty\}.
\end{array}
$$
Indeed, for each $\phi\in S(\Lambda,\mu)$, its  integral $\hat J_\phi$
with respect to the spectral measure $\hat{P}$ and the multiplication
operator $Q_\phi$ coincide, that is,
$Q_\phi=\hat J_\phi=\int \phi\,d \hat{P}$,
and then we can translate every result of the theory of  integration
with respect to spectral measures to the context of multiplication
operators on direct integrals.

Moreover, each measurable family of operators 
$\{T(\l):\H_\l\to\H_\l\}_{\l\in \Lambda}$ defines an 
operator on ${\H}_{\mu,N}$, say $T$, given by
$Tf:= \int_\Lambda^\oplus T(\l)f(\l)\,d\mu(\l)$, ($f=f(\l)\in D(T)$).
The operators of this form are called  {\it
decomposable operators}.
In particular, when $T(\l)=\phi(\l)I_\l$, with $I_\l$ the identity operator on $\H_\l$, then
$T$ is just the multiplication operator $Q_\phi$.
A decomposable operator $T$ is self-adjoint 
(unitary, normal, orthogonal projection) if and only if the operators
$T(\l)$ are self-adjoint (unitary, normal, orthogonal projection) on
$\H_\l$ for $\mu$-a.e. $\l\in\Lambda$.

Like for spectral measures $P$ the spectral
type $[P]$ and the multiplicity function $N_P$ are the unitary
invariants, the measure type $[\mu]$ and the  function of
dimension $N$ are the unitary invariants of direct integrals
$\H_{\mu,N}$, i.e., they determine the
corresponding structure uniquely up to unitary equivalence. 
Now, given a spectral measure space $(\Lambda,\A,\H,P)$ and
a direct integral $\H_{\mu,N}=\int_\Lambda^\oplus \H_\l\,d\mu(\l)$  
defined on $(\Lambda,\A)$, 
a {\it structural isomorphism} among both structures is a  unitary
operator $V$ from $\H$ onto $\H_{\mu,N}$ verifying 
\begin{equation}\label{19xx}
V\, J_\phi= Q_\phi\, V,\quad \phi\in S(\Lambda,P)=S(\Lambda,\mu).
\end{equation}
In particular, $V P ( E )= \hat{P}( E )V$, for each $ E \in{\A}$.
The conditions under which a spectral measure space and a  direct
integral are structurally isomorphic and the explicit form of a
structural isomorphism are given in the following \cite[Th.7.4.1]{BS}:

\begin{thm}\label{pA2}

A spectral measure space  $(\Lambda,\A,\H,P)$ and a direct 
integral $\H_{\mu,N}=\int_\Lambda^\oplus \H_\l\,d\mu(\l)$ defined on
$(\Lambda,\A)$ are structurally isomorphic if and only if
$$
[P]=[\mu]\quad {\rm and}\quad N_P=N\,\,[\mu]\mbox{-}{\rm a.e.}
$$
In such case, a structural isomorphism $V:\H\to\H_{\mu,N}$ among them 
is defined by 
\begin{equation}\label{19}
V^{-1} h=\bigoplus_{j=1}^m \left(\int_{\Lambda}
\left(e_j(\l),\sqrt{\frac{d\mu}{d\mu_{g_1}}(\l)}\;h(\l)
\right)_\l \,dP(\l)\right)\,g_j \,,
\end{equation}
where $\{e_j(\l)\}_{j=1}^m$ is a measurable orthonormal basis in
$\H_{\mu,N}$ and $\{g_k\}_{k=1}^m$ is a generating system of  $\H$ with
respect to $P$ such that:

\smallskip
i.) $[P]=[\mu_{g_1}]\succ [\mu_{g_2}]\succ\dots$; 

\smallskip
ii.) if $1\le j\le k\le m$ ($m=[P]\mbox{-}{\rm ess\, sup\,} N_P$), then the 
measure
$\mu_{g_k}$ coincides with $\mu_{g_k}$ on its support $\Lambda(g_k)$,
which means that $\mu_{g_j}|_{\Lambda(g_k)}=\mu_{g_k}$.
\end{thm}

The functional version of von Neumann spectral theorem  \cite{VN49}
is a direct consequence of the above facts. 

Finally, recall that, given a generating system $\{g_j\}_{j=1}^m$ of $\H$ with respect to a spectral measure $P$, the space $\H$ decomposes into the orthogonal sum (\ref{odh}), so that each $h\in \H$ can be written as $h=\oplus_{j=1}^m h_j$, where $h_j\in \H_{g_j}$. This fact together with property {\it (vii)} of Proposition \ref{pA1} imply that for each $h\in\H$ there exists a unique family of functions $\tilde{h}_j\in L^2(\Lambda,\mu_{g_j})$  such that 
\begin{equation}\label{erkk34}
h=\bigoplus_{j=1}^m h_j = \bigoplus_{j=1}^m J_{\tilde{h}_j}g_j=
\bigoplus_{j=1}^m \big[\int_{\Lambda} \tilde{h}_j\,dP\big]\,g_j.
\end{equation}

Theorem \ref{tA4} below collects some results proved in \cite{GG}. The first one describes the $\l j$-components of $h$ and $Vh$ in terms of certain Radom-Nikodym derivatives associated with the generating system. The second one shows that it is possible to obtain a spectral decomposition  ``a la
Dirac" by using direct integrals of Hilbert spaces. 

\begin{thm}\label{tA4}
Under the conditions of Theorem \ref{pA2} we have: 

\smallskip
(i) For each
$h\in\H$ and $j\in\{1,2,\ldots,m\}$, the following identities are satisfied $\mu$-a.e.:
\begin{equation}\label{28}
\begin{array}{l}
\ds 
\sqrt{\frac{d\mu_{g_j}}{d\mu}(\l)}\,\tilde{h}_j(\l)= \big(e_j(\l),[Vh](\l)\big)_{\l} =
\\[3ex]
\hspace{0ex} \ds
= \sqrt{\frac{d\mu}{d\mu_{g_1}}(\l)}\,\frac{d\mu_{g_j,h}}{d\mu}(\l)
= \sqrt{\frac{d\mu_{g_1}}{d\mu}(\l)}\,\frac{d\mu_{g_j,h}}{d\mu_{g_1}}(\l)
= \sqrt{\frac{d\mu_{g_j}}{d\mu}(\l)}\,\frac{d\mu_{g_j,h}}{d\mu_{g_j}}(\l)\,.
\end{array}
\end{equation}

\smallskip
(ii) For each $f,h\in\H$ and $E\in\A$, 
\begin{equation}
\big(f,P(E)h\big)_\H=\sum_{j=1}^m \int_E \big([Vf](\l), e_j(\l)\big)_\l\;
\big(e_j(\l), [Vh](\l)\big)_\l\,d\mu(\l) \label{20}
\end{equation}
\end{thm}


\subsection{Locally Convex Equipments of Spectral Measures}\label{app3}

By a {\it topological vector space} ({\bf tvs}) we mean a pair ($\Phi,\tau_\Phi$),
 where $\Phi$ is a vector space over the complex field ${\mathbb C}$  and $\tau_\Phi$ a locally convex linear topology on $\Phi$. In what follows $\Phi^\t$ denotes the space of continuous {\it antilinear} mappings (functionals) from $\Phi$ into $\mathbb C$.
The space $\Phi^\t$ is also a complex vector space and we can endow it with its own topology \cite{RR}.

For the (anti)dual pair $(\Phi,\Phi^\t)$ we denote the action of  $F^\t\in\Phi^\t$
on $\phi\in\Phi$ as $\<\phi|F^\t\>$. This bracket is linear to the right and antilinear to
the left, just like the scalar product of Hilbert spaces. As usual, we write $\<F^\t|\phi\>=\<\phi|F^\t\>^\ast$, where $\ast$ denotes the complex conjugate.

In order to approach the question of existence of generalized
eigenvectors  for a normal operator defined on a Hilbert space
$\H$ and with continuous spectrum,  one can propose to rig the
Hilbert space $\H$ with a pair
$(\Phi,\Phi^\t)$ such that
$\Phi$ is a dense subspace of $\H$ and $\Phi^\t$ contains a complete set of
generalized eigenvectors for the operator in question or, better, for the associated spectral measure.  
To be precise:
\begin{defn}\rm
We say that a tvs $(\Phi,\tau_\Phi)$ {\it rigs or equips the
spectral measure space} $(\Lambda,{\A},\H,P)$ when the following conditions hold:

\smallskip
i.) There exists a one to one linear mapping $I:\Phi\longmapsto \H$ with
range dense  in $\H$. Identifying
each $\phi\in\Phi$ with its image $I(\phi)$, we can assume that
$\Phi\subset \H$ is a dense subspace of $\H$ and $I$ the canonical
injection from $\Phi$ into $\H$.

\smallskip
ii.) There exists a $\sigma$-finite measure $\mu$ on
$(\Lambda,\A)$, a set $\Lambda_0\subset \Lambda$ with zero $\mu$
measure and a family of vectors in $\Phi^\t$ of the form
\begin{equation}\label{36}
\Big\{{\l k}^\t\in\Phi^\t \,:\, \l\in \Lambda\backslash \Lambda_0,\,
k\in\{1,2,\ldots,m\}\Big\}, 
\end{equation}
where $m\in\{\infty,1,2,\ldots\}$, such that
\begin{equation}\label{37}
(\phi, P ( E )\varphi)_{\H}=\sum_{k=1}^m \int_ E   \<\phi|{\l
k}^\t\>\, \<{\l k}^\t|\varphi\>\,d\mu(\l),\quad
(\phi,\varphi\in\Phi,\,\,E \in{\A}).
\end{equation}
In particular, if $E=\Lambda$, then, $P(E)=I_\H$, the identity on $\H$, and
$$
(\phi,\varphi)_{\H}=\sum_{k=1}^m \int_\Lambda  \<\phi|{\l
k}^\t\>\, \<{\l k}^\t|\varphi\>\,d\mu(\l),\quad
(\phi,\varphi\in\Phi).
$$ 

A family of the form (\ref{36}) satisfying (\ref{37})
is called {\it a complete system of generalized
eigenvectors} (also called {\it Dirac kets}) of the spectral measure space $(\Lambda,{\A},\H,P)$. 
\end{defn}

Each direct integral $\H_{\mu,N}$ associated to
$(\Lambda,{\A},{\H}, P )$ as in Theorem \ref{pA2} along with
one of its measurable orthonormal bases $\{e_k(\l)\}_{k=1}^{N(\l)}$ or,
equivalently, each generating system $\{g_k\}_{k=1}^m$ in ${\H}$ with
respect to $P$ verifying conditions i.) and ii.) of Theorem
\ref{pA2}, provide a rigging $(\Phi,\tau_\Phi,\mu,\{\l k^\t\})$. This
rigging is characterized by the following properties \cite{GG}:

\begin{itemize}
\item[(i)] 
The subspace $\Phi$ is dense in $\H$ and is given by 

$$
\Phi=\{\phi\in {\H}: {\rm exists}\,\,
\frac{d\mu_{\phi,g_k}}{d\mu_{g_k}}(\l)<\infty,\,\,\forall \l\in
\Lambda\backslash \Lambda_0,\,\,\forall k\in\{1,2,\ldots,N(\l)\}\},
$$
where $\Lambda_0$ is a subset of $\Lambda$ with $\mu$ zero measure (or
equivalently, $P$ zero measure).

\item[(ii)] 
The complete family of antilinear functionals on $\Phi$ fulfilling
(\ref{37}) is of the form
$$
\Big\{\l k^\t\,:\, \l\in\Lambda\backslash \Lambda_0,\,
k\in\{1,2,\ldots,N(\l)\}\Big\},
$$
where the action of $\l k^\t$ over each $\phi\in\Phi$ is given by:
\begin{equation}\label{39}
\< \phi|\l k^\t\> = \big(e_k(\l),[Vh](\l)\big)_{\l} = 
\sqrt{\frac{d\mu_{g_k}}{d\mu}(\l)}\,
\frac{d\mu_{g_k,h}}{d\mu_{g_k}}(\l)\,.
\end{equation}

\item[(iii)]
The topological antidual space $\Phi^\t$ is the vector space spanned by the set $\{\l k^\t\}$. The topology $\tau_\Phi$ is the weak topology $\sigma(\Phi,\Phi^\t)$, i.e. the coarsest one compatible with the dual pair $(\Phi,\Phi^\t)$.  
The topology $\tau_\Phi$ is generated by the family of seminorms
$$
\phi\mapsto |\<\phi|\l k^\t\>|, \quad  \l\in \Lambda\backslash 
\Lambda_0,\,\,k\in\{1,2,\ldots,N(\l)\}.
$$
\end{itemize}

This type of riggings is {\it minimal} 
in the sense that no topology on $\Phi$ coarser than that given in (iii) --except for
the indeterminacy derived by the zero $\mu$ measure set
$\Lambda_0$--  can rig the spectral measure space $(\Lambda,{\A},{\H},
P )$.
In the opposite side we find the so-called {\it universal} equipments, those with Hilbert-Schmidt inductive and nuclear topologies which, due to the inductive and nuclear versions of the spectral  theorem, rig every ``regular" spectral measure space, {\it c.f.} \cite{GG1}.


\section{Transformation Theory}\label{sctf}

Let us pass to deal with the main subject of this work: the transformation formulas between two given representations described by spectral measures $(\Lambda,{\A},\H,P)$ and
$(\Xi,{\B},\H,Q)$, here called $P$ and $Q$ for brevity. The
elements we are going to handle are the following:

\begin{itemize}

\item[({\bf A})]
The scalar measures derived of the spectral measures $P$ and $Q$,
$$
\mu_{f,g}^P(E):=(f,P(E)g)_\H
\quad {\rm and}\quad
\mu_{f,g}^Q(F):=(f,Q(F)g)_\H,
$$
where $f,g\in\H$, $E\in{\A}$ and $F\in{\B}$.
Recall that when $f=g$ we write $\mu_f^P:=\mu_{f,f}^P$ and
$\mu_f^Q:=\mu_{f,f}^Q$. (See Section \ref{app1}.)

\item[({\bf B})]
Generating systems $\{g_j\}_{j=1}^m$ and $\{g'_k\}_{k=1}^{m'}$ 
of $\H$ with respect to $P$ and $Q$, respectively, which satisfy the
following properties (see Section \ref{app1} and Theorem \ref{pA2}):
\begin{itemize}
\item[(a)] 
$[ P ]=[\mu^P_{g_1}]\succ [\mu^P_{g_2}]\succ\cdots$ and $[ Q
]=[\mu^Q_{g'_1}]\succ [\mu^Q_{g'_2}]\succ\cdots$.
\item[(b)] 
If $1\leq j\leq k\leq m$, then $\mu^P_{g_j}|_{\Lambda(g_k)}=\mu^P_{g_k}$ and, if $1\leq j\leq
k\leq m'$, then
$\mu^Q_{g'_j}|_{\Lambda(g'_k)}=\mu^Q_{g'_k}$.
\end{itemize}

\item[({\bf C})]
Direct integrals of Hilbert spaces (see Section \ref{app2})
$$
\H_{\mu,N}=\int_\Lambda \H_\l\,d\mu(\l)\quad {\rm and} \quad
\H_{\mu',N'}=\int_\Xi \H_\xi\,d\mu'(\xi)
$$
associated, respectively, to $P$ and $Q$ as in Theorem \ref{pA2} 
through orthonormal measurable bases in both
direct integrals
$$
\{e_j(\l)\}_{j=1}^m\subset\H_{\mu,N} 
\quad{\rm and }\quad 
\{e'_k(\xi)\}_{k=1}^{m'}\subset\H_{\mu',N'}
$$
and corresponding structural isomorphisms 
$$
{\V}:\H\to\H_{\mu,N}
\quad{\rm and }\quad 
{\V}':\H\to\H_{\mu',N'}\,.
$$
We shall suppose without loss of generality that
$$
[\mu]=[P]=[\mu_g^P]
\quad{\rm and}\quad
[\mu']=[Q]=[\mu^Q_{g_1}].
$$

\item[({\bf D})]
A locally convex topological vector space ({\bf tvs}) $[\Phi,\tau_\Phi]$ that rigs both
spectral  measures and two complete systems of generalized eigenvectors
$$
\{\l j^\t\}_{(\l,j)\in\Lambda\t\{1,2,\ldots,m\}}\subset\Phi^\t \quad {\rm
and}\quad 
\{\xi k^{\t}\}_{(\xi,k)\in\Xi\t\{1,2,\ldots,m'\}}\subset\Phi^\t
$$
for both spectral measures $P$ and $Q$, being the decompositions given with respect to the
scalar  measures $\mu$ and $\mu'$, respectively. 
(See Section \ref{app3}.)
\end{itemize}

Identifying generalized eigenvectors with the elements of the orthonormal measurable bases given in ({\bf C}), transformation formulas 
(\ref{2g}) in direct integrals should be similar, for each $h\in\H$, to

\begin{equation}
\label{D.19.5444id}
\big(e'_k(\xi),[{\V}'h](\xi)\big)_{\xi}= \sum_{j=1}^m
\int_\Lambda 
\< e'_k(\xi)|e_j(\l)\>\,
\big(e_j(\l),[{\V}h](\l)\big)_{\l}\,d\mu(\l)\,.
\end{equation}
On the other hand, when we consider the locally convex  equipment
$[\Phi,\tau_\Phi]$ given in ({\bf D}), transformation formulas
should be similar, for each $\psi\in\Phi$, to
\begin{equation}
\label{D.19.5444rhs}
\< \xi k^\t|\psi\>= \sum_{j=1}^m \int_\Lambda 
\< \xi k^\t|\l j^\t\>\, \< \l j^\t|\psi\>\, d\mu(\l)\,. 
\end{equation} 
But the integrals kernels of these formulas, written as 
$$
\< e'_k(\xi)|e_j(\l)\>
\quad {\rm and}\quad
\< \xi k^\t|\l j^\t\>\,,
$$
in general, do not have algebraic sense.
Indeed, $\< e'_k(\xi)|e_j(\l)\>$ combines elements belonging to
different spaces: $e'_k(\xi)\in\H_\xi$ and
$e_j(\l)\in\H_\l$. And, though $\xi k^\t$ and $\l j^\t$ are both in $\Phi^\t$, in general neither
$\xi k^\t$ nor $\l j^\t$ belong to $\Phi$ and therefore $\< \xi k^\t|\l
j^\t\>$ is meaningless in the dual pair $(\Phi,\Phi^\t)$.

We shall analyze equations (\ref{D.19.5444id}) and (\ref{D.19.5444rhs}) 
at the same time because the fo\-llo\-wing relations are satisfied: 
\begin{equation}\label{emm0P}
\< \phi|\l j^\t\>= \big([{\V}\phi](\l),e_j(\l)\big)_{\H_\l}=
\sqrt{\frac{d\mu^P_{g_k}}{d\mu}(\l)}\,\frac{d\mu^P_{\phi,g_k}}{d\mu^P_{g_k}}
(\l),
\quad (\phi\in\Phi),
\end{equation}

\begin{equation}\label{emm0Q}
\< \phi|\xi k^\t\>= \big([{\V'}\phi](\xi),e'_k(\xi)\big)_{\H_\xi}=
\sqrt{\frac{d\mu^Q_{g'_k}}{d\mu'}(\xi)}\,\frac{d\mu^Q_{\phi,g'_k}}{d\mu^Q_{g'_k}}
(\xi),
\quad (\phi\in\Phi).
\end{equation}
(See Equation (\ref{39}).)
Theorems \ref{p18} and \ref{p16} below mean a first approach to the problem. Both results consider  Vitali approach to Radon-Nikodym derivatives (see Appendix), so that the following Lemma is all we need to prove them.

\begin{lem}\label{l17}
Let $\H_{\mu,N}$ and $\H_{\mu',N'}$ be two direct integrals corresponding, respectively, to the spectral measure spaces $(\Lambda,{\A},\H,P)$ and $(\Xi,{\B},\H,Q)$
as in ({\bf C}). Then for all $F\in{\B}$ and each pair of elements
$f,h\in\H$ we have the following identities:
\begin{equation}\label{tle44.nc}
\begin{array}{rcl}
\ds \mu^Q_{f,h}(F) 
& = & \ds 
\sum_{j=1}^m \int_\Lambda 
\big( [{\V}Q(F)f](\l),e_j(\l) \big)_{\l}\,
\big(e_j(\l),[{\V}h](\l)\big)_{\l}\,d\mu(\l) \\[2ex]
& = & \ds 
\sum_{j=1}^m \int_\Lambda 
\big( [{\V}f](\l),e_j(\l) \big)_{\l}\,
\big(e_j(\l),[{\V}Q(F)h](\l)\big)_{\l}\,d\mu(\l).
\end{array}
\end{equation}
\end{lem}

\begin{proof}
Let $f$ and $h$ be two elements of $\H$ whose
decompositions (\ref{erkk34}) are given by $f=\oplus_j f_j$ and
$h=\oplus_j h_j$, where
$f_j,h_j\in\H_{g_j}$. The elements of the family
$\{h_{j}\}_{j\in\{1,2,\ldots,m\}}$ are spectrally orthogonal with respect to
$P$ and $\sum_j||h_{j}||^2<\infty$. Thus we have  
$$
\mu^P_{h}(E)=\mu^P_{\oplus_j h_{j}}( E )=\sum_j\mu^P_{h_{j}}( E ), \quad
(E\in{\A}),
$$
and the same is valid for $\{f_{j}\}_{j\in\{1,2,\ldots,m\}}$.
Therefore,
\begin{equation}\label{mcmp44}
\mu^P_{f,h}(E)=\sum_j\mu^P_{f,h_{j}}( E )=
\sum_j\mu^P_{f_j,h}( E )=\sum_j\mu^P_{f_j,h_{j}}( E ),
\quad (E\in{\A}).
\end{equation}
As in (\ref{erkk34}), for each $j\in\{1,2,\ldots,m\}$, denote by $\tilde{f}_{j}$
the function of $L^2(\Lambda,\mu^P_{g_j})$ such that 
\begin{equation}\label{mcmp44x}
f_{j}=J_{\tilde{f}_{j}}g_j=\big[\int_\Lambda \tilde{f}_{j}\, dP\big]g_j.
\end{equation}
From (\ref{mcmp44}) and (\ref{mcmp44x}), being $\V$ unitary, it follows  
\begin{equation}\label{ettdcdeg}
\begin{array}{rcl}
\mu^Q_{f,h}(F) & = & \ds \sum_{j=1}^m \big( f_j,Q(F)h \big)_\H
=  \sum_{j=1}^m \big( J_{\tilde{f}_j}g_j,Q(F)h \big)_\H
\\[2ex] 
& = & \ds \sum_{j=1}^m \big( g_j,J_{\tilde{f}_j^\ast}Q(F)h \big)_\H
=  \sum_{j=1}^m \big( {\V}g_j,{\V}J_{\tilde{f}_j^\ast}Q(F)h
\big)_{\H_{\mu,N}}
\end{array}
\end{equation}
Now, in the light of (\ref{19}), 
$\ds \V g_j=\sqrt{\frac{d\mu^P_{g_1}}{d\mu}(\l)}\,e_j(\l)$.
This fact together with (\ref{19xx}) imply  
$$
\begin{array}{l}
\ds\sum_{j=1}^m \big( {\V}g_j,{\V}J_{\tilde{f}_j^\ast}Q(F)h
\big)_{\H_{\mu,N}} = 
\\
\ds =\sum_{j=1}^m \int_\Lambda
\sqrt{\frac{d\mu^P_{g_1}}{d\mu}(\l)}\,
\tilde{f}_j^\ast(\l)\,
\big( e_j(\l), [\V Q(F)h](\l)\big)_{\l}\,d\mu(\l).
\end{array}
$$
Finally, from (\ref{28}),
$$
\tilde{f}_j^\ast(\l)=
\sqrt{\frac{d\mu}{d\mu^P_{g_1}}(\l)}\,\big([{\V}f_j](\l),e_j(\l)\big)_{\l}=
\sqrt{\frac{d\mu}{d\mu^P_{g_1}}(\l)}\,\big([{\V}f](\l),e_j(\l)\big)_{\l}
,\quad \mu\mbox{-}{\rm c.s.} 
$$
Substituting in (\ref{ettdcdeg}) we obtain the
first  identity of (\ref{tle44.nc}). The second one is obtained
interchanging the roles of $f$ and $h$ in (\ref{ettdcdeg}).
\end{proof}

\begin{thm}\label{p18}
Let $\H_{\mu,N}$ and $\H_{\mu',N'}$ be two direct integrals corresponding, res\-pec\-ti\-ve\-ly, to the spectral measure spaces $(\Lambda,{\A},\H,P)$ and $(\Xi,{\B},\H,Q)$
as in ({\bf C}). Su\-ppo\-se the measure $\mu'$ has a Vitali system and let
$\{F_n\}_{n=1}^\infty$ be a sequence of sets of ${\B}$ admitting a
contraction to a point $\xi\in\Xi$ and such that $\mu'(F_n)\neq0$ ($n\in\mathbb N$). Then, for all $h\in\H$ and 
$k\in\{1,2,\ldots,m'\}$, 
\begin{equation}\label{ettdc4idns.nc}
\begin{array}{l}
\ds \big([{\V}' h](\xi),e'_k(\xi)\big)_\xi =\\[1ex] 
 =  \ds \sqrt{\frac{d\mu'}{d\mu^Q_{g'_k}}(\xi)}\, 
\lim_{n\to\infty} 
\sum_{j=1}^m \int_{\Lambda} 
\big(e_j(\l),[{\V}g'_{k}](\l)\big)_{\l}\,
{\big( [{\V}Q(F_n)h](\l),e_j(\l) \big)_{\l}\over \mu'(F_n)}\,
d\mu(\l) \\[3ex]
 =  \ds 
\sqrt{\frac{d\mu'}{d\mu^Q_{g'_k}}(\xi)}\, 
\lim_{n\to\infty} 
\sum_{j=1}^m \int_\Lambda 
{\big(e_j(\l),[{\V}Q(F_n)g'_{k}](\l)\big)_{\l}\over \mu'(F_n)} \,
\big( [{\V}h](\l),e_j(\l) \big)_{\l}\,
d\mu(\l).
\end{array}
\end{equation}
\end{thm}

\begin{proof}
By Equation (\ref{39}), save for a set $\Lambda_0$ of zero 
$\mu'$-measure, we have 
\begin{equation}\label{drnsvit4.0}
\big([{\V}' h](\xi),e'_k(\xi)\big)_\xi=
\sqrt{\frac{d\mu^Q_{g'_k}}{d\mu'}(\xi)}\,
\frac{d\mu^Q_{h,g'_k}}{d\mu^Q_{g'_k}}(\xi)=
\sqrt{\frac{d\mu'}{d\mu^Q_{g'_k}}(\xi)}\,
\frac{d\mu^Q_{h,g'_k}}{d\mu'}(\xi).
\end{equation}
By hypothesis the measure $\mu'$ have a Vitali system. 
Since for each pair $f,h\in\H$ the measure $\mu^Q_{f,h}$ is absolutely 
continuous with respect to $\mu'$, every Vitali system for $\mu'$ is
also a Vitali system for
$\mu^Q_{f,h}$.  Thus, the Vitali-Lebesgue theorem (Theorem \ref{tA7} of Appendix) implies that
for  
$\mu'$-a.a. $\xi\in\Xi$,
if $F_1,F_2,\ldots$ is a sequence of measurable sets of ${\B}$ that 
admits a contraction to  $\xi$, then
\begin{equation}\label{drnsvit4}
\frac{d\mu^Q_{h,g'_k}}{d\mu'}(\xi)=\lim_{n\to\infty}
\frac{\mu^Q_{h,g'_k}(F_n)}{\mu'(F_n)}.
\end{equation}
Substituting $\mu^Q_{h,g'_k}(F_n)$ in (\ref{drnsvit4}) by the corresponding right hand sides of (\ref{tle44.nc}) and putting them into (\ref{drnsvit4.0}) we get formulae (\ref{ettdc4idns.nc}). 
\end{proof}

Appealing to equalities (\ref{emm0P}), (\ref{emm0Q}) and 
(\ref{28}), Theorem \ref{p18} can be translated directly into the framework
of locally convex equipments. The only additional fact one needs is that functionals belong to $\Phi^\t$. This happens, for example, when $\tau_\Phi$ is finer than the subspace topology induced by $\H$ on $\Phi$, so that $\Phi\subseteq\H\subseteq\Phi^\t$. 

\begin{thm}\label{p16}
Consider a locally convex tvs
$[\Phi,\tau_\Phi]$ rigging both spectral measure spaces $(\Lambda,{\A},\H,P)$ and $(\Xi,{\B},\H,Q)$ as in ({\bf D}) such that $\tau_\Phi$ is finer than the subspace topology induced by $\H$ on $\Phi$. Assume the measure $\mu'$ has a Vitali system and let
$\{F_n\}_{n=1}^\infty$ be a sequence of sets in ${\B}$ admitting a
contraction to the point $\xi\in\Xi$ and such that $\mu'(F_n)\neq0$ ($n\in\mathbb N$). Then, for all $\phi\in\Phi$ and $k\in\{1,2,\ldots,m'\}$,
\begin{equation}\label{ettdc4evtns}
\ds \<\phi | \xi k^\t \>  = 
\sqrt{\frac{d\mu'}{d\mu^Q_{g'_k}}(\xi)}\, 
\lim_{n\to\infty} 
\sum_{j=1}^m \int_\Lambda 
\sqrt{\frac{d\mu_{g_k}}{d\mu}(\l)}\,{[\widetilde{Q(F_n)g'_{k}}](\l)\over \mu'(F_n)}
\,\< \phi | \l j^\t\>\,d\mu(\l).
\end{equation}
In particular, when each $g'_{kj}$ belongs to $\Phi$ and $Q(F_n)\Phi\subseteq\Phi$ for all $n\in\mathbb N$,
\begin{equation}\label{ettdc4evtns.1}
\begin{array}{rcl}
\ds \<\phi | \xi k^\t \> 
& = & \ds \sqrt{\frac{d\mu'}{d\mu^Q_{g'_k}}(\xi)}\, 
\lim_{n\to\infty} \sum_{j=1}^m \int_{\Lambda} 
\< \l j^\t | g'_{k}\> \,\frac{\< Q(F_n)\phi | \l j^\t\>}{\mu'(F_n)}\,d\mu(\l) \\[3ex]
& = & \ds
\sqrt{\frac{d\mu'}{d\mu^Q_{g'_k}}(\xi)}\, 
\lim_{n\to\infty} \sum_{j=1}^m \int_\Lambda 
\frac{\< \l j^\t | Q(F_n)g'_{k}\>}{\mu'(F_n)}\,\< \phi | \l j^\t\>\,d\mu(\l).
\end{array}
\end{equation}
\end{thm}

\begin{rem}\rm
Under some additional assumptions formulae simplify notably. For example,
let us take $\mu=\mu_{g_1}$ and $\mu'=\mu_{g'_1}$, so that normalization factors disappear, and suppose the spectral measures $P$ and $Q$ commute.
\begin{itemize}
\item
If, moreover, $P$ is simple, the generating system with respect to $P$ has only one element $g$ and the set of operators commuting with every projection $P(E)$,
($E\in{\A}$),  is just the algebra of multiplication operators with
respect to $P$. In particular, for each projection
$Q(F)$, ($F\in{\B}$), there exists a function $\gamma_{F}\in
L^\infty(\Lambda,\mu^P_g)$ such that $Q(F)=J_{\gamma_{F}}$.
Moreover, being $Q(F)$ an orthogonal projection, $\gamma_{F}$ must be the
characteristic  function of certain set $E_{F}\in{\A}$, i.e.
$\gamma_{F}=\chi_{E_{F}}$. Therefore, in this case (\ref{ettdc4evtns}) has the following form:
\begin{equation}\label{ettdc4ft}
\< \phi|\xi k^{\t}\> 
=  \lim_{n\to\infty} 
\int_{E_{F_n}} \frac{\tilde{g}'_k(\l)}{\mu'(F_n)}\,\<\phi|\l^\t\>\,d\mu(\l),
\end{equation}
In particular, if $g'_k\in\Phi$,  
\begin{equation}\label{ettdc4fft}
\< \phi|\xi k^{\t}\> =
 \lim_{n\to\infty} 
\int_{E_{F_n}} \frac{\<\l^\t|g'_k\>}{\mu'(F_n)}\,\<\phi|\l^\t\>\,d\mu(\l).
\end{equation}
Equation (\ref{ettdc4idns.nc}) simplifies in a similar way:
\begin{equation}\label{ettdc4ft.dih}
\big([{\V}' h](\xi),e'_k(\xi)\big)_\xi
=  \lim_{n\to\infty} 
\int_{E_{F_n}} \frac{\tilde{g}'_k(\l)}{\mu'(F_n)}\,\big( [{\V}h](\l),e(\l) \big)_{\l}\,d\mu(\l),
\end{equation}
where $e(\l)$ is the unique element of the orthogonal measurable basis in $\H_{\mu,N}$.
\item
For arbitrary spectral measure $P$, not necessarily simple, in
general, the ope\-ra\-tors $Q(F)$ are not multiplication but decomposable operators of the form
$$
Q(F)=\int_\Lambda^\oplus [Q(F)](\l)\,d\mu(\l),
$$
where $[Q(F)](\l)$ is an orthogonal projection in $\H_\l$ for $\mu$-a.e.
$\l\in\Lambda$. In this case (\ref{ettdc4idns.nc}) takes the form:
\begin{equation}\label{ettdc4idns}
\begin{array}{l}
\ds \big([{\V}' h](\xi),e'_k(\xi)\big)_\xi =\\[1ex] 
 =  \ds \lim_{n\to\infty} \sum_{j=1}^m \int_{\Lambda} 
\big(e_j(\l),[{\V}g'_{k}](\l)\big)_{\l}\,
\frac{\big( [Q(F_n)](\l)[{\V}h](\l),e_j(\l) \big)_{\l}}{\mu'(F_n)}\,d\mu(\l) 
\\[2ex]
 \ds = \lim_{n\to\infty} \sum_{j=1}^m \int_\Lambda 
\frac{\big(e_j(\l),[Q(F_n)](\l)[{\V}g'_{k}](\l)\big)_{\l}}{\mu'(F_n)}\,
\big( [{\V}h](\l),e_j(\l) \big)_{\l}\,d\mu(\l).
\end{array}
\end{equation}
\end{itemize}
\end{rem}

 When, as in Theorem \ref{p18}, one considers  direct integrals
of Hilbert spaces, $\H\simeq\H_{\mu,N}\simeq\H_{\mu',N'}$, being Hilbert spaces, one can identify them as usual with their topological antidual spaces $\H^\t\simeq\H^\t_{\mu,N}\simeq\H^\t_{\mu',N'}$. The 
isomorphism $\H\simeq\H^\t$ is linear
whereas the isomorphism between $\H$ and its topological dual $\H'$ is antilinear --this is the reason why antiduals instead duals are considered--.

Fix $(\xi,k)\in\Xi\t\{1,2,\ldots,m'\}$ and for each pair $(n,l)$, with $n,l\in\mathbb N$ and $1\leq l<m$, define the functional $\tilde\Gamma^\t_{nl}:\H \to {\mathbb C}$ by
\begin{equation}\label{fun.ant.Gammat.deg.nc}
\tilde\Gamma^\t_{nl}(h):= 
 \sqrt{\frac{d\mu'}{d\mu^Q_{g'_k}}(\xi)}\,\sum_{j=1}^l
\int_{\Lambda}  
{\big(e_j(\l),[{\V}Q(F_n)g'_{k}](\l)\big)_{\l}\over\mu'(F_n)}
 \,\big( [{\V}h](\l),e_j(\l)\big)_\l\,d\mu(\l)\,,
\end{equation}
($h\in\H$). Really $\tilde\Gamma^\t_{nl}$ depends on $(\xi,k)$ also, but we drop this fact from the notation.
Obviously $\{\tilde\Gamma^\t_{nl}\}\subset\H^\t$ and (\ref{ettdc4idns.nc}) is equivalent to
the simple or  pointwise convergence on $\H$ of the sequence $\{\tilde\Gamma^\t_{nl}\}$
to the functional $\big([{\V}'\cdot](\xi),e'_k(\xi)\big)_\xi$ or, in other words, its
convergence on $\H^\t$ with respect to the weak topology $\sigma(\H^\t,\H)$ --the usual ``weak limit", say $\wlim$, in Hilbert spaces--, i.e.
$$
\big([{\V}'\cdot](\xi),e'_k(\xi)\big)_\xi=\wlim_{(n,l)\to(\infty,m)} \tilde\Gamma^\t_{nl}.
$$

In a similar way, under the hypothesis of Theorem \ref{p16}, the functionals 
$\Gamma^\t_{nl}:\Phi \to {\mathbb C}$ defined by
\begin{equation}\label{fun.ant.Gammaxx}
\Gamma^\t_{nl}(\phi):= \sqrt{\frac{d\mu'}{d\mu^Q_{g'_k}}(\xi)}\, \sum_{j=1}^l
\int_\Lambda 
\sqrt{\frac{d\mu_{g_k}}{d\mu}(\l)}\,{[\widetilde{Q(F_n)g'_{k}}](\l)\over
\mu'(F_n)}\, \< \phi|\l j^\t\>\,
d\mu(\l)\,,
\end{equation}
($\phi\in\Phi$), belong to $\Phi^\t$, for each $n\in\mathbb N$ and $1\leq l<m$, and (\ref{ettdc4evtns}) says that $\{\Gamma^\t_{nl}\}$ converges to $\xi k^\t$ on $\Phi^\t$ with respect to the weak topology $\sigma(\Phi^\t,\Phi)$, that is,
$$
\xi k^\t=\lim_{(n,l)\to(\infty,m)} \Gamma^\t_{nl}\quad
{\rm in\ \ } [\Phi^\t,\sigma(\Phi^\t,\Phi)]\,.
$$

Clearly the sequence $\{\tilde\Gamma^\t_{nl}\}$ cannot converge to $\big([{\V}'\cdot](\xi),e'_k(\xi)\big)_\xi$ in strong sense on $\H^\t\simeq\H$. Indeed, even in the simplest case described by 
(\ref{ettdc4ft.dih}), where the index $l$ disappears and the sequence $\{\tilde\Gamma^\t_{n}\}$ may be identified with the sequence of elements of $\H_{\mu,N}$ given by
$$
\left\{{\chi_{E_{F_n}}(\l)\,\tilde{g}'_k(\l) \over \mu'(F_n)}\,e(\l)\right\}_{n\in\mathbb N},
$$
since $\ds \lim_{n\to\infty} \mu'(F_n)=\mu'(\xi)=0$ and the characteristic functions $\chi_{E_{F_n}}$ take values $0$ and $1$ only, if  $\{\tilde\Gamma^\t_{n}\}$ converges in strong sense, its limit must be the zero function, but this contradicts 
 the fact that $e'_k(\xi)$ belongs to an orthonormal
measurable basis of $\H_{\mu',N'}$. 

In what follows we establish conditions for the convergence of $\{\Gamma^\t_{nl}\}$ to $\xi k^\t$ with respect to topologies finer than weak topology on $\Phi^\t$.


\subsection{$\bf \tau_{pc}(\Phi^\t,\Phi)$-convergence.}\label{ss.2.1.1}

Among the most important results on barrelled tvs we mention the
Banach-Steinhaus theorem \cite{RR}, which involves the topology
$\tau_\sigma$ of simple or pointwise convergence and the topology
$\tau_{pc}$ of uniform convergence on precompact sets 
of $\Phi$ defined on the  space ${\mathcal L}^\t(\Phi,\Psi)$ of all
$(\tau_\Phi,\tau_\Psi)$-continuous
(anti)linear mappings from $\Phi$ into $\Psi$, where $[\Psi,\tau_\Psi]$
is any other locally convex tvs.

\begin{thm}\label{t3}
{\bf [Banach-Steinhaus]}
Let $\Phi$ and $\Psi$ be
locally  convex tvs, $\Phi$ in addition barrelled. If $(T_\alpha)$ is a net
in
${\mathcal L}^\t(\Phi,\Psi)$ which is
$\tau_\sigma$-bounded\footnote{
If $\Phi$ is a barrelled 
tvs then the $\sigma(\Phi^\t,\Phi)$-bounded sets and the
$\beta(\Phi^\t,\Phi)$-bounded sets coincide.
}
and which converges pointwise to some $T\in\Psi^\Phi$, then
$T\in{\mathcal L}^\t(\Phi,\Psi)$ and $(T_\alpha)$ converges to $T$ with respect to
the topology $\tau_{pc}$.
\end{thm}

In the study of the transformation theory on locally convex equipments
the  relevant part of the  Banach-Steinhaus theorem is the
$\tau_{pc}$-convergence of the net, because the other part of the result is
implicit in the construction of the equipment.

\begin{cor}\label{p5}
Let $[\Phi,\tau_\Phi]$ be a locally convex tvs 
 rigging both spectral measures as in ({\bf D}) and
such that $[\Phi,\tau_\Phi]$ is a barrelled space whose topology is finer
than that induced by $\H$ on
$\Phi$. Then, under the conditions of Theorem \ref{p16}, for all
$k\in\{1,2,\ldots,m'\}$ and for $\mu'$-a.a. $\xi\in\Xi$, we have  
\begin{equation}\label{l.tpc.cd}
\xi k^{\t} 
= \lim_{(n,l)\to(\infty,m)} \Gamma^\t_{nl}
\quad
{\rm in\ \ } [\Phi^\t,\tau_{pc}(\Phi^\t,\Phi)]\,.
\end{equation}
\end{cor}

\begin{proof}
Equation (\ref{ettdc4evtns}) assures the sequence
$\{\Gamma^\t_{nl}\}\subset \Phi^\t$ is $\sigma(\Phi^\t,\Phi)$-bounded and
converges pointwise to $\xi k^\t$.
By Banach-Steinhaus theorem, $\{\Gamma^\t_{nl}\}$ converges
to $\xi k^\t$ in $[\Phi^\t,\tau_{pc}(\Phi^\t,\Phi)]$.
\end{proof}

\subsection{$\bf \beta(\Phi^\t,\Phi)$-convergence.}\label{ss.2.1.3}

Now we pass to study under what conditions we can substitute the 
$\tau_{pc}$-convergence of the sequence $\{\Gamma^\t_{nl}\}$ by the convergence with respect to the
strong topology $\beta(\Phi^\t,\Phi)$. We focus attention on properties of barrelled
and semi-reflexive tvs.

It is well known that a locally convex tvs $[\Phi,\tau_\Phi]$ is 
semi-reflexive if and only if the Mackey topology $\mu(\Phi^\t,\Phi)$
and the strong topology $\beta(\Phi^\t,\Phi)$ coincide in $\Phi^\t$
\cite{RR}. Thus, for a semi-reflexive tvs $\Phi$ we have
\begin{equation}\label{ig.top.etsr}
\mu(\Phi^\t,\Phi)=\tau_{pc}(\Phi^\t,\Phi)=\beta(\Phi^\t,\Phi).
\end{equation}

The following result is immediate from Corollary \ref{p5}
and identities (\ref{ig.top.etsr}):

\begin{cor}\label{p9}
Let $[\Phi,\tau_\Phi]$ be a barrelled and
semi-reflexive tvs\footnote{A barrelled and semi-reflexive locally convex tvs $\Phi$ admit the
following  characterizations \cite{RR}: 
(a) $\Phi$ is reflexive;
(b) $\Phi$ is semi-reflexive and quasi-barrelled;
(c) $\Phi$ is quasi-barrelled and weakly quasi-complete.
}
that rigs both spectral measures as in ({\bf
D})  with a topology finer than that induced by $\H$ on
$\Phi$. Then, under the conditions of Theorem \ref{p16}, for all $k\in\{1,2,\ldots,m'\}$ and for $\mu'$-a.a. $\xi\in\Xi$, we have 
$$
\xi k^{\t} 
= \lim_{(n,l)\to(\infty,m)} \Gamma^\t_{nl}
\quad
{\rm in\ \ } [\Phi^\t,\beta(\Phi^\t,\Phi)]\,.
$$
\end{cor}

In particular, the conclusions of Proposition \ref{p9} are satisfied when $[\Phi,\tau_\Phi]$ is a Montel space, a Fr\'echet-Montel space\footnote{Since a Fr\'echet space $\Phi$ is
a Montel  space if and only if it is separable and every
$\sigma(\Phi',\Phi)$-convergent sequence in $\Phi'$ is
$\beta(\Phi',\Phi)$-convergent, and formula (\ref{ettdc4ft}) implies
the weak convergence of the sequence of functionals $\Gamma_n$,
for Fr\'echet-Montel spaces we can deduce directly the conclusions of Corollary \ref{p9} from Theorem \ref{p16}. } 
or a locally convex tvs such that the strong antidual space $\Phi_\beta^\t$
is nuclear\footnote{
If $\Phi$ is a locally convex tvs such that the strong antidual space $\Phi_\beta^\t$
is nuclear, then $\Phi$ is 
semi-Montel and quasi-barrelled, that is, $\Phi$ is a Montel space
\cite{PI72}.
}. 


\section{Final Remarks}\label{sect.32}

Transformation theory has been introduced in Section \ref{sctf} on the basis of Vitali approach to Radon-Nikodym derivatives, in terms of which generalized eigenfunctions are given. Any other approach to the concept of Radon-Nikodym derivative, as for example martingale theory or generalized Cauchy-Stieltjes and Poisson integrals, must lead to similar results.

For example, let us take $(\Xi, {\B},\mu')$ a probability space  and let $\nu$ be a finite measure on $(\Xi, {\B})$. Assume there is a sequence
$\Pi_n$, $n=1,2,\dots$, of partitions of $\Xi$ into a finite number of pairwise disjoint measurable subsets of positive measure $\mu'$, such that each partition $\Pi_{n+1}$  is finer than
$\Pi_n$ and the $\sigma$-algebra generated by $\bigcup_{n=1}^{\infty}\Pi_n$ coincides with ${\B}$. 
Put 
$$ 
X_n(\l)=\sum_{F\in\Pi_n}{\nu(F)\over \mu'(F)}\chi_F(\l)\,,
$$ 
and denote by ${\B}_n$  the sub-$\sigma$-algebra of ${\B}$ generated by
$\Pi_n$. Then $\{(X_n,{\B}_n): n=1,2,\ldots\}$ is a martingale that converges $\mu'$-a.e. to an integrable limit $X$, which coincides with the usual Radon-Nikodym derivative
$d\nu/d\mu'$ provided $\nu$ is absolutely continuous with respect to $\mu'$ \cite{Billin}. 
Clearly the functionals $^M\tilde\Gamma^\t_{nl}$ and $^M\Gamma^\t_{nl}$ given by 
$$
\begin{array}{rccl}
^M\tilde\Gamma^\t_{nl}: & \H & \longrightarrow & {\mathbb C} 
\\[2ex]
 & h & \mapsto &
\ds \sqrt{\frac{d\mu'}{d\mu^Q_{g'_k}}(\xi)}\,\sum_{j=1}^l
\sum_{F\in\Pi_n}{\chi_F(\xi)}\t
\\[1ex]
&&&\ds\t\int_{\Lambda}  
{\big(e_j(\l),[{\V}Q(F)g'_{k}](\l)\big)_{\l}\over\mu'(F)}
 \,\big( [{\V}h](\l),e_j(\l)\big)_\l\,d\mu(\l)
\end{array}
$$
and
$$
\begin{array}{rccl}
^M\Gamma^\t_{nl}: & \Phi & \longrightarrow & {\mathbb C} \\[2ex]
 & \phi & \mapsto &
\ds \sqrt{\frac{d\mu'}{d\mu^Q_{g'_k}}(\xi)}\, \sum_{j=1}^l
\sum_{F\in\Pi_n}{\chi_F(\xi)}\t
\\[1ex]
&&&\ds\t\int_\Lambda 
\sqrt{\frac{d\mu_{g_k}}{d\mu}(\l)}\,{[\widetilde{Q(F)g'_{k}}](\l)\over
\mu'(F)}\, \< \phi|\l j^\t\>\,
d\mu(\l)
\end{array}
$$
have the same properties than the functionals $\tilde\Gamma^\t_{nl}$ and $\Gamma^\t_{nl}$ defined in Equations (\ref{fun.ant.Gammat.deg.nc}) and (\ref{fun.ant.Gammaxx}), respectively.
Approximation by finite $\sigma$-algebras, as in this exam\-ple, can be characterized in terms of  separability of $L^1(\mu')$ \cite{MALL}.

\subsection*{Acknowledgments}
We thank to Profs. I.E. Antoniou, A. Bohm, M. Castagnino, F.
L\'opez Fdez-Asenjo, M. N\'u\~nez and Z. Suchanecki for useful discussions. 
This work was supported by JCyL-project VA013C05 (Castilla y Le\'on) and MEC-project FIS2005-03989 (Spain).


\appendix

\section{Vitali Systems.}

Let $(\Lambda,{\A},\mu)$ be a measure space such that
for each $\l\in\Lambda$ the set $\{\l\}$ is measurable and
$\mu(\l)=0$. 

\begin{defn}\rm
A {\it Vitali system} for $(\Lambda,{\A},\mu)$ is a family of measurable sets ${\mathcal V}\subseteq \A$ such that:
\begin{itemize}
\item[(i)]
Being given a measurable set $E\in\A$ and $\varepsilon >0$, there exist
a countable family of Vitali sets $A_1,A_2,\dots$ such that
$$
E \subset \bigcup_{n=1}^\infty A_n\,, \hskip0.7cm  \mu\left(
\bigcup_{n=1}^\infty A_n \right) <\mu(E)+\varepsilon
$$ 
\item[(ii)]
Each $A\in\A$  has a {\it border}, i.e., a zero $\mu$-measure
set $\partial A$ such that
\begin{itemize}
\item[(a)]
If $\l\in A\setminus(A\cap\partial A)$, then any Vitali set
containing $\l$ with measure sufficiently small is contained in
$A\setminus(A\cap\partial A)$.
\item[(b)]
If $x\notin A\cup\partial A$, then any Vitali set
containing $\l$ with measure sufficiently small has no common point
with
$A\cup\partial A$.
\end{itemize}

\item[(iii)]
Let $E\subset \Lambda$ be a set admitting a covering by Vitali sets
$\B\subset V$ such that for each $\l\in E$ and each $\varepsilon>0$,
there exists
a set $A_\varepsilon(x)\in\B$ with $\mu[A_\varepsilon(x)]<\varepsilon$ and $\l\in
A_\varepsilon(x)$. Then, $E$ can be covered, to within a set of $\mu$ zero
measure, by countably many disjoint sets $A_j \in \B$.
\end{itemize}
\end{defn}

\begin{ex}\rm
The Lebesgue measure on $({\mathbb R}^n,{\B})$, where $\B$
is the Borel $\sigma$-algebra, has as particular Vitali system the set of
all closed cubes.
\end{ex}

\begin{defn}\rm
A sequence of measurable sets $\{E_1,E_2,\dots\}\subset\A$ {\it admits a contraction}  to a point $\l_0\in\Lambda$  when:
\begin{itemize}
\item[(i)]
For each $E_n$ in the sequence, there is a Vitali set $A_n$ such
that $\l_0\in A_n$ and $\lim_{n\to\infty} \mu(A_n)=0$.
\item[(ii)]
There exists a positive constant $c$ ($c>0$) such that
$$
\mu(E_n)\ge c\mu(A_n)\,,\hskip0.7cm (n\in \mathbb N)
$$
\end{itemize}
\end{defn}

\begin{defn}\rm
Let $\nu$ be a countably additive function from $\A$ into ${\mathbb R}$. 
The {\it de\-ri\-va\-ti\-ve of $\nu$ at the point $\l_0$  with
respect to the Vitali system $\mathcal V$} is given by
$$
D_{\mathcal V}(\l_0) =\lim_{\varepsilon\to
0}\frac{\nu[A_\varepsilon(\l_0)]}{\mu[A_\varepsilon(\l_0)]}
$$
(provided that the limit exists), where $A_\varepsilon(\l_0)$ is any Vitali
set with $\mu$ measure smaller than $\varepsilon$ containing $\l_0$.
\end{defn}

The basic result on differentiation with respect to a Vitali system is due to Vitali and Lebesgue \cite[Th.10.1]{SG}:

\begin{thm}\label{tA7}
{\bf [Vitali-Lebesgue]}
Let $(\Lambda,{\A},\mu)$ be
a measure space, $\mathcal V$ a Vitali system on $\Lambda$ and $\nu$ a
mapping from
$\A$ into $\mathbb R$ with the property of being absolutely continuous with
respect to
$\mu$. Then, the derivative of $\nu$ with respect to $\mathcal V$ exists,
save for a set of $\mu$ zero measure. This is given at each point
$\l_0$ by

$$
D_{\mathcal V}(\l_0)=\lim_{n\to \infty} \frac{\nu(E_n)}{\mu(E_n)}
$$
where $E_1,E_2,\dots$ is a sequence of measurable sets admitting a
contraction to $\l_0$. This derivative coincides with the
Radon-Nikodym derivative of $\nu$ with respect to $\mu$,
$d\nu/d\mu$.
\end{thm}



\end{document}